\documentclass[12pt,twoside]{amsart}

\usepackage{graphicx}
\newtheorem{theorem}{\bf Theorem}

\newtheorem{definition}{\bf Definition}

\begin{document}
\vspace{1.3cm}
\title{On existence and asymptotic behaviour of solutions of a fractional integral equation with linear
modification of the argument} \maketitle
\begin{center}
\author{{\large Mohamed Abdalla Darwish\footnote{Permanent address:Department of Mathematics, Faculty of
Science, Alexandria University at Damanhour, 22511 Damanhour, Egypt. E-mail:darwishma@yahoo.com.}}\\
Department of Mathematics, Massachusetts Institute of Technology\\77 Massachusetts Ave., Cambridge,
MA 02139-4307, USA.\\
E-mail: darwish@math.mit.edu}

\end{center}
\begin{abstract}
We study the solvability of a quadratic integral equation of fractional order with
linear modification of the argument. This equation is considered in the Banach space
of real functions defined, bounded and continuous on an unbounded interval. Moreover, we will
obtain some asymptotic characterization of solutions.

\vskip 0.4 true cm \noindent AMS Mathematics Subject
Classification : 45G10, 45M99, 47H09.

\end{abstract}

\pagestyle{myheadings} \markboth{\centerline {M.A. Darwish}}
{\centerline{On asymptotic behaviour of solutions of a fractional integral equation}}
\section{Introduction}
The study of differential and integral equations with
a modified argument is relatively new, it was initiated only
in the past thirty years or so. These equations arise in the modeling of problems from the
natural and social sciences such as biology, physics and economics, see
\cite{{BaMi}, {BeDa5}, {CaLoSa}, {CaDy}, {Dal}, {Du}, {Ku}, {Me},
{Mu03}, {Mu99}} and references therein. On the other hand, The first serious attempt to give a logical definition of a
fractional derivative is due to Liouville, see \cite{Hil} and
references therein. Today differential and integral equations of fractional order play a very
important role in describing numerous events and problems of the real world. For example, many problems in
mechanics, physics, economics and other fields led to differential and integral equations of fractional order
(cf. \cite{{BaRz2}, {BaOr}, {BaOl1}, {DMB}, {DDSA}, {DCAA}, {Da}, {Hil}, {Po}, {Sam}}).

In this paper, we will study the quadratic integral equation of fractional order with linear
modification of the argument in the integral, namely
\begin{equation}
\label{e1}x(t)=a(t)+\frac{f(t,x(t))}{\Gamma(\alpha)}\;\int_0^t
\frac{u(t,s,x(s),x(\lambda s))}{(t-s)^{1-\alpha}}\;
ds,
\end{equation}
where $t\in\mathbb{R}_+$ and $0<\alpha,\;\lambda<1$. Throughout $a:\mathbb{R}_+\rightarrow\mathbb{R}$,
$f:\mathbb{R}_+\times\mathbb{R}\rightarrow\mathbb{R}$ and
$u:\mathbb{R}_+\times\mathbb{R}_+\times\mathbb{R}\times\mathbb{R}\rightarrow\mathbb{R}$ are functions which satisfy
special assumptions, see Section\ref{sec3}. Let us recall that the
function $f=f(t,x)$ involves in Eq.(\ref{e1}) generates the
superposition operator $F$ defined by
\begin{equation}
\label{sup}(Fx)(t)=f(t,x(t)),
\end{equation}
where $x=x(t)$ is an arbitrary function defined on $\mathbb{R}_+$, see \cite{ApZa}.

In the case $\alpha=1$, $a(t)=1$, $f(t,x)=x$ and
$u(t,s,x,y)=\frac{t}{t+s}\;\phi(s)\;x$, where $\phi$ is a continuous function and $\phi(0)=0$, Eq.(\ref{e1})
has the form
\begin{equation}
\label{18} x(t)=1+x(t)\int_0^t\frac{t}{t+s}\;\phi(s)\;x(s)\;ds.
\end{equation}
Eq.(\ref{18}) is the Volterra counterpart of the famous quadratic integral equation of
Chandrasekhar type considered in many papers and monographs (cf.
\cite{{Ar}, {BaLeWa}, {Ch}, {HuKhZh}} for instance). Some Problems
considered in the theory of radiative transfer, in the theory of
neutron transport and in the kinetic theory of gases lead to
quadratic integral equations (cf. \cite{{BaLeWa}, {BaRz1}, {Ch},
{DJQT}, {De}, {HuKhZh}, {Ke}, {Le}}).

The aim of this paper is to prove the existence of solutions of Eq.(\ref{e1}) in the space of real
functions, defined, continuous and bounded on an unbounded interval. Moreover, we will
obtain some asymptotic characterization of solutions of Eq.(\ref{e1}). Our proof depends on suitable
combination of the technique of measures
of noncompactness and the Schauder fixed point principle.

It is worthwhile mentioning that up to now the work of J. Bana{\'s} and D. O'Regan \cite {BaOr} is the only paper concerning
with the study of quadratic integral equation of factional order in the space of real functions defined, continuous and
bounded on an unbounded interval.
\section{Notation and auxiliary facts}
\setcounter{equation}{0}
This section is devoted to collecting some definitions and results
which will be needed further on.
First we recall the definition of the Riemann-Liouville fractional integral \cite{{Po}, {Sam}}.
\begin{definition}{\rm
Let $f\in L_1(a,b)$, $0\leq a<b<\infty$, and let $\alpha>0$
be a real number. The Riemann-Liouville fractional integral of order $\alpha$ of the function $f(t)$ is
defined by
$$I^\alpha f(t)=\frac{1}{\Gamma(\alpha)}
\int_0^t\frac{f(s)}{(t-s)^{1-\alpha}}\;ds,\;a<t<b.$$}
\end{definition}

Now, let $(E,\|.\|)$ be an infinite dimensional Banach space with zero element $0$. Let $B(x,r)$ denotes
the closed ball centered at $x$ with radius $r$. The symbol
$B_r$ stands for the ball $B(0,r)$.

If $X$ is a subset of $E$, then $\overline{X}$ and
$Conv X$ denote the closure and convex closure of $X$,
respectively. Moreover, we denote by ${\mathcal{M}}_E$ the family
of all nonempty and bounded subsets of $E$ and by
${\mathcal{N}}_E$ its subfamily consisting of all relatively
compact subsets.

Next we give the concept of a measure of noncompactness \cite{BaGo}:
\begin{definition}\label{d1}{\rm
A mapping $\mu:{\mathcal{M}}_E\rightarrow\mathbb{R}_+=[0,\infty)$ is
said to be a measure of noncompactness in $E$ if it
satisfies the following conditions:
\begin{enumerate}
\item[$1)$] The family ${\rm ker}\mu=\{X\in{\mathcal{M}}_E:\mu(X)=0\}$ is nonempty and ${\rm
Ker}\mu\subset{\mathcal{N}}_E$. \item[$2)$] $X\subset
Y\Rightarrow \mu(X)\leq\mu(Y)$. \item[$3)$] $\mu(\overline{X})=\mu(Conv
X)=\mu(X)$. \item[$4)$] $\mu(\lambda
X+(1-\lambda)Y)\leq\lambda\;\mu(X)+(1-\lambda)\;\mu(Y)$ for
$0\leq\lambda\leq 1$. \item[$5)$] If $X_n\in{\mathcal{M}}_E$,
$X_n=\overline{X}_n$, $X_{n+1}\subset X_n$ for $n=1,\;2,\;3,\;...$ and
$\lim\limits_{n\rightarrow\infty}\mu(X_n)=0$ then $\cap_{n=1}^\infty
X_n\neq\phi$.
\end{enumerate}}
\end{definition}
The family ${\rm ker}\mu$ described above is called the kernel of
the measure of noncompactness $\mu$. Let us observe that the intersection set $X_\infty$ from $5)$
belongs to ${\rm ker}\mu$. In fact, since $\mu(X_\infty)\leq\mu(X_n)$ for every $n$ then we have that
$\mu(X_\infty)=0$.

In what follows we will work in the Banach space $BC(\mathbb{R}_+)$
consisting of all real functions defined, bounded and continuous on $\mathbb{R}_+$. This space is
equipped with the standard norm
$$\|\;x\;\|=\sup\{|x(t)|:t\geq 0\}.$$

Now, we recollect the construction of the measure of noncompactness
in $BC(\mathbb{R}_+)$ which will be used in the next section
(see \cite{{Ba}, {BaOl}}).

Let us fix a nonempty and bounded subset $X$ of $BC(\mathbb{R}_+)$ and numbers $\varepsilon>0$ and
$T>0$. For arbitrary function $x\in X$ let us denoted by $\omega^T(x,\varepsilon)$ the
modulus of continuity of the function $x$ on the interval $[0,T]$, i.e.
$$\omega^T(x,\varepsilon)={\rm sup}\{|x(t)-x(s)|:t,\;s
\in[0,T],\;|t-s|\leq\varepsilon\}$$ Further, let us put
$$\omega^T(X,\varepsilon)={\rm sup}\{\omega^T(x,\varepsilon):x\in X\},$$
$$\omega_0^T(X)=\lim_{\varepsilon\rightarrow
0}\omega^T(X,\varepsilon)$$
and
$$\omega_0^\infty(X)=\lim_{T\rightarrow
\infty}\omega_0^T(X,\varepsilon).$$
Moreover, for a fixed number $t\in\mathbb{R}_+$ let us defined
$$X(t)=\{x(t):x\in X\}$$
and
$${\rm diam}X(t)=\sup\{|x(t)-y(t)|:x,\;y\in X\}.$$
Finally, let us define the function $\mu$ on the family
${\mathcal{M}}_{BC(\mathbb{R}_+)}$ by
\begin{equation}
\label{f1}
\mu(X)=\omega_0^\infty(X)+c(X),
\end{equation}
where $c(X)=\limsup\limits_{t\rightarrow\infty}{\rm diam}X(t).$
The function $\mu$ is a measure of noncompactness in the space $BC(\mathbb{R}_+)$, see \cite{{Ba}, {BaOl1}}.
Let us mention that the kernel ${\rm ker}\mu$ of the measure $\mu$ consists of all sets
$X\in{\mathcal{M}}_{BC(\mathbb{R}_+)}$ such
that functions from $X$ are locally equicontinuous on $\mathbb{R}_+$ and vanish uniformly at infinity, i.e.
for any $\varepsilon>0$ there exists $T>0$ such that $|x(t)|\leq\varepsilon$ for all $x\in X$ and for any $t\geq T$.
This property will permit us to characterize solutions of Eq.(\ref{e1}).
\section{Existence Theorem}\label{sec3}
\setcounter{equation}{0}
In this section we will study Eq.(\ref{e1}) assuming that the
following hypotheses are satisfied:
\begin{enumerate}
\item[$(h_1)$] $a:\mathbb{R}_+\rightarrow\mathbb{R}$ is a continuous and bounded function on $\mathbb{R}_+$.
\item[$(h_2)$] $f:\mathbb{R}_+\times\mathbb{R}\rightarrow\mathbb{R}$ is continuous and there exists a continuous function
$m(t)=m:\mathbb{R}_+\rightarrow\mathbb{R}_+$ such that
$$|f(t,x)-f(t,y)|\leq m(t)|x-y|$$
for all $x,\;y\in\mathbb{R}$ and for any $t\in\mathbb{R}_+$.
\item[$(h_3)$]
$u:\mathbb{R}_+\times\mathbb{R}_+\times{\mathbb{R}}\times{\mathbb{R}}\rightarrow{\mathbb{R}}$
is a continuous function. Moreover, there exist a function $n(t)=n:\mathbb{R}_+\rightarrow\mathbb{R}_+$ being continuous
on $\mathbb{R}_+$ and a function $\Phi:\mathbb{R}_+\times\mathbb{R}_+\rightarrow\mathbb{R}_+$ being
continuous and nondecreasing on $\mathbb{R}_+$ with $\Phi(0,0)=0$ and such
that
$$|u(t,s,x_2,y_2)-u(t,s,x_1,y_1)|\leq n(t)\Phi(|x_2-x_1|,|y_2-y_1|)$$ for all $t,\;s\in\mathbb{R}_+$ such
that $t\geq s$ and for all
$x_i,\;y_i\in\mathbb{R}$ $(i=1,\;2)$.

For further purpose let us define the function $u^*:\mathbb{R}_+\rightarrow\mathbb{R}_+$
by $u^*(t)=\max\{|u(t,s,0,0)|:0\leq s\leq t\}$.
\item[$(h_4)$] The functions $\phi,\;\psi,\;\xi,\;\eta:\mathbb{R}_+\rightarrow\mathbb{R}_+$ defined by
$\phi(t)=m(t)\;n(t)\;t^\alpha$, $\psi(t)=m(t)\;u^*(t)\;t^\alpha$, $\xi=n(t)|f(t,0)|t^\alpha$ and
$\eta(t)=u^*(t)|f(t,0)|t^\alpha$ are bounded on $\mathbb{R}_+$ and the functions $\phi$ and $\xi$ vanish at infinity, i.e.
$\lim\limits_{t\rightarrow\infty}\phi(t)=\lim\limits_{t\rightarrow\infty}\xi(t)=0$.
\item[$(h_5)$] There
exists a positive solution $r_0$ of the inequality
\begin{equation}
\label{inq}\|a\|\Gamma(\alpha+1)+[\phi^*r\Phi(r,r)+\psi^*r+\xi^*\Phi(r,r)+\eta^*]\leq
r\;\Gamma(\alpha+1)
\end{equation}
 and $\phi^*\Phi(r_0,r_0)+\psi^*<\Gamma(\alpha+1)$, where
 $\phi^*=\sup\{\phi(t):t\in\mathbb{R}_+\}$, $\psi^*=\sup\{\psi(t):t\in\mathbb{R}_+\}$,
 $\xi^*=\sup\{\xi(t):t\in\mathbb{R}_+\}$ and $\eta^*=\sup\{\eta(t):t\in\mathbb{R}_+\}$.
\end{enumerate}
Now, we are in a position to state and prove our main result.
\begin{theorem}\label{mth}
{\rm Let the hypotheses $(h_1)-(h_5)$ be satisfied.
Then Eq.(\ref{e1})
has at least one solution $x\in BC(\mathbb{R}_+)$ such that $x(t)\rightarrow 0$
as $t\rightarrow\infty$.}
\end{theorem}
{\bf Proof:} Denote by $\mathcal{F}$ the operator associated with
the right-hand side of equation (\ref{e1}), i.e., equation (\ref{e1}) takes the form
$$x={\mathcal{F}}x,$$
where
\begin{equation}
\label{c3} ({\mathcal{F}}\;x)(t)=a(t)+(Fx)(t)\cdot({\mathcal{U}}\;x)(t),
\end{equation}
and
\begin{equation}
\label{c33} ({\mathcal{U}}\;x)(t)=\frac{1}{\Gamma(\alpha)}\int_0^t
\frac{u(t,s,x(s),x(\lambda s))}{(t-s)^{1-\alpha}}\;ds.
\end{equation}
Solving Eq.(\ref{e1}) is equivalent to finding a fixed
point of the operator $\mathcal{F}$ defined on the space $BC(\mathbb{R}_+)$.

We claim that for any function $x\in BC(\mathbb{R}_+)$ the operator $\mathcal{F}$ is continuous on $\mathbb{R}_+$.
To establish this claim it suffices to show that if
$x\in BC(\mathbb{R}_+)$ then ${\mathcal{U}}x$ is continuous function on $\mathbb{R}_+$, thanks $(h_1)$ and $(h_2)$.
For, take an arbitrary $X\in BC(\mathbb{R}_+)$ and fix $\varepsilon>0$ and $T>0$. Assume that
$t_1,\;t_2\in\mathbb{R}_+$ are such that $|t_2-t_1|\leq\varepsilon$. Without loss of generality we can assume that
$t_2>t_1$. Then we get
{\small\begin{eqnarray} \nonumber
|({\mathcal{U}}x)(t_2)-({\mathcal{U}}x)(t_1)|&=&\left|\frac{1}{\Gamma(\alpha)}\int_0^{t_2}
\frac{u(t_2,s,x(s),x(\lambda s))}{(t_2-s)^{1-\alpha}}\;ds\right.\\
\nonumber\;&\;&\;\;\;\;\;\;\;-\left.\frac{1}{\Gamma(\alpha)}\int_0^{t_1}
\frac{u(t_1,s,x(s),x(\lambda s))}{(t_1-s)^{1-\alpha}}\;ds\right|\\
\nonumber\;&\leq&\left|\frac{1}{\Gamma(\alpha)}\int_0^{t_2}
\frac{u(t_2,s,x(s),x(\lambda s))}{(t_2-s)^{1-\alpha}}\;ds\right.\\
\nonumber\;&\;&\;\;\;\;\;\;\;-\left.\frac{1}{\Gamma(\alpha)}\int_0^{t_1}
\frac{u(t_2,s,x(s),x(\lambda s))}{(t_2-s)^{1-\alpha}}\;ds\right|\\
\nonumber\;&\;&\;\;+\left|\frac{1}{\Gamma(\alpha)}\int_0^{t_1}
\frac{u(t_2,s,x(s),x(\lambda s))}{(t_2-s)^{1-\alpha}}\;ds\right.\\
\nonumber\;&\;&\;\;\;\;\;\;\;-\left.\frac{1}{\Gamma(\alpha)}\int_0^{t_1}
\frac{u(t_1,s,x(s),x(\lambda s))}{(t_2-s)^{1-\alpha}}\;ds\right|\\
\nonumber\;&\;&\;\;+\left|\frac{1}{\Gamma(\alpha)}\int_0^{t_1}
\frac{u(t_1,s,x(s),x(\lambda s))}{(t_2-s)^{1-\alpha}}\;ds\right.\\
\nonumber\;&\;&\;\;\;\;\;\;\;-\left.\frac{1}{\Gamma(\alpha)}\int_0^{t_1}
\frac{u(t_1,s,x(s),x(\lambda s))}{(t_1-s)^{1-\alpha}}\;ds\right|\\
\nonumber\;&\leq&\frac{1}{\Gamma(\alpha)}\int_{t_1}^{t_2}\frac{|u(t_2,s,x(s),x(\lambda s))|}
{(t_2-s)^{1-\alpha}}\;ds\\
\nonumber\;&\;&\;\;+\frac{1}{\Gamma(\alpha)}\int_0^{t_1}
\frac{|u(t_2,s,x(s),x(\lambda s))-u(t_1,s,x(s),x(\lambda s))|}{(t_2-s)^{1-\alpha}}\;ds\\
\nonumber\;&\;&+\frac{1}{\Gamma(\alpha)}\int_0^{t_1}
|u(t_1,s,x(s),x(\lambda(s)))|[(t_1-s)^{\alpha-1}-(t_2-s)^{\alpha-1}]ds.
\end{eqnarray}}
Therefore, if
\begin{eqnarray}
\nonumber
\omega_d^T(u,\varepsilon)=\sup\{|u(t_2,s,y_1,y_2)-u(t_1,s,y_1,y_2)|:\;s,\;t_1,\;t_2\in[0,T],
\;\;\;\;\;\;\;\;\;\;\;\;\;\;\;\\
\nonumber\;\;\;\;\;\;\;\;\;\;\;\;\;\;\;\;\;\;\;\;\;\;\;\;\;\;\;\;\;\;\;t_1
\geq s,\;t_2\geq s,\;|t_2-t_1|
\leq\varepsilon,\;{\rm and}\;y_1,\;y_2\in[-d,d]\}
\end{eqnarray}
then we obtain
{\small\begin{eqnarray} \nonumber
|({\mathcal{U}}x)(t_2)-({\mathcal{U}}x)(t_1)|\;&\leq&
\frac{1}{\Gamma(\alpha)}\int_{t_1}^{t_2}\frac{|u(t_2,s,x(s),x(\lambda s))-u(t_2,s,0,0)|+|u(t_2,s,0,0)|}
{(t_2-s)^{1-\alpha}}\;ds\\
\nonumber\;&\;&\;\;+\frac{1}{\Gamma(\alpha)}\int_0^{t_1}
\frac{\omega_{\|x\|}^T(u,\varepsilon)}{(t_2-s)^{1-\alpha}}\;ds\\
\nonumber\;&\;&+\frac{1}{\Gamma(\alpha)}\int_0^{t_1}
[|u(t_1,s,x(s),x(\lambda s))-u(t_1,s,0,0)|+|u(t_1,s,0,0)|]\\
\nonumber\;&\;&\;\;\;\;\;\;\;\;\;\;\;\;\;\;\;\;\;\;\;\;\;\;\;\times[(t_1-s)^{\alpha-1}-(t_2-s)^{\alpha-1}]ds\\
\nonumber\;&\leq&
\frac{1}{\Gamma(\alpha)}\int_{t_1}^{t_2}\frac{n(t_2)\Phi(|x(s)|,|x(\lambda s)|)+u^*(t_2)}
{(t_2-s)^{1-\alpha}}\;ds\\
\nonumber\;&\;&\;\;+\frac{\omega_{\|x\|}^T(u,\varepsilon)}{\Gamma(\alpha+1)}[t_2^\alpha-(t_2-t_1)^\alpha]\\
\nonumber\;&\;&+\frac{1}{\Gamma(\alpha)}\int_0^{t_1}
[n(t_1)\Phi(|x(s)|,|x(\lambda s)|)+u^*(t_1)]\\
\nonumber\;&\;&\;\;\;\;\;\;\;\;\;\;\;\;\;\;\;\;\;\;\;\;\;\;\;\times[(t_1-s)^{\alpha-1}-(t_2-s)^{\alpha-1}]ds\\
\nonumber\;&\leq&
\frac{n(t_2)\Phi(\|x\|,\|x\|)+u^*(t_2)}{\Gamma(\alpha+1)}(t_2-t_1)^\alpha
+\frac{\omega_{\|x\|}^T(u,\varepsilon)}{\Gamma(\alpha+1)}t_1^\alpha\\
\nonumber\;&\;&+\frac{n(t_1)\Phi(\|x\|,\|x\|)+u^*(t_1)}{\Gamma(\alpha+1)}[t_1^\alpha-t_2^\alpha+(t_2-t_1)^\alpha].
\end{eqnarray}}
Thus
\begin{equation}\label{es}
\omega^T({\mathcal{U}}x,\varepsilon)\leq\frac{1}{\Gamma(\alpha+1)}\left\{
2\varepsilon^\alpha[\hat{n}(T)\Phi(\|x\|,\|x\|)+\hat{u}(T)]+T^\alpha\omega_{\|x\|}^T(u,\varepsilon)\right\},
\end{equation}
where we denoted
$$\hat{n}(T)=\max\{n(t):t\in[0,T]\}$$
and
$$\hat{u}(T)=\max\{u^*(t):t\in[0,T]\}.$$

In view of the uniform continuity of the function $u$ on $[0,T]\times[0,T]\times[-\|x\|,\|x\|]\times[-\|x\|,\|x\|]$ we have that
$\omega_{\|x\|}^T(u,\varepsilon)\rightarrow 0$ as $\varepsilon\rightarrow
0$. From the above inequality we infer that the function ${\mathcal{U}}x$ is continuous on the interval $[0,T]$ for
any $T>0$. This yields the continuality of ${\mathcal{U}}x$ on $\mathbb{R}_+$, and
consequently, the function ${\mathcal{F}}x$ is continuous on $\mathbb{R}_+$.

Now, we show that ${\mathcal{F}}x$ is bounded on $\mathbb{R}_+$. Indeed, in view of our hypotheses for arbitrary
$x\in BC(\mathbb{R}_+)$ and for a fixed
$t\in\mathbb{R}_+$ we have
{\small\begin{eqnarray} \nonumber
|({\mathcal{F}}x)(t)|&\leq&\left|a(t)+\frac{f(t,x(t))}{\Gamma(\alpha)}
\int_0^{t}\frac{u(t,s,x(s),x(\lambda s))}{(t-s)^{1-\alpha}}\;ds\right|\\
\nonumber\;&\leq&\|a\|+\frac{1}{\Gamma(\alpha)}[|f(t,x(t))-f(t,0)|+|f(t,0)|]\\
\nonumber&\;&\;\;\;\;\;\;\times\int_0^{t}
\frac{|u(t,s,x(s),x(\lambda s))-u(t,s,0,0)|+|u(t,s,0,0)|}{(t-s)^{1-\alpha}}\;ds\\
\nonumber\;&\leq&\|a\|+\frac{m(t)\|x\|+|f(t,0)|}{\Gamma(\alpha)}\int_0^{t}
\frac{n(t)\Phi(|x(s)|,|x(\lambda s)|)+u^*(t)}{(t-s)^{1-\alpha}}\;ds\\
\nonumber\;&\leq&\|a\|+\frac{m(t)\|x\|+|f(t,0)|}{\Gamma(\alpha+1)}\;[n(t)\;\Phi(\|x\|,\|x\|)+u^*(t)]\;t^\alpha\\
\nonumber\;&=&\|a\|+\frac{1}{\Gamma(\alpha+1)}[\phi(t)\|x\|\Phi(\|x\|,\|x\|)+\psi(t)\|x\|+\xi(t)\Phi(\|x\|,\|x\|)+\eta(t)].
\end{eqnarray}}
Hence, ${\mathcal{F}}x$ is bounded on $\mathbb{R}_+$, thanks hypothesis $(h_4)$. This assertion in conjunction with
the continuity of ${\mathcal{F}}x$ on $\mathbb{R}_+$ allows us to conclude that the operator $\mathcal{F}$ maps
$BC(\mathbb{R}_+)$ into itself. Moreover, from the last estimate we have
$$\|{\mathcal{F}}x\|\leq\|f\|+\frac{1}{\Gamma(\alpha+1)}\;[\phi^*\|x\|\Phi(\|x\|,\|x\|)+\psi^*\|x\|
+\xi^*\Phi(\|x\|,\|x\|)+\eta^*].$$
Linking this estimate with hypothesis $(h_5)$ we deduce that there exists $r_0>0$ such that the operator
$\mathcal{F}$ transforms the ball $B_{r_0}$ into itself.

Now, we prove that the operator $\mathcal{F}$ is continuous on
the ball $B_{r_0}$. To do this, let us fix $\varepsilon>0$ and take $x,\;y\in B_{r_0}$ such that $\|x-y\|\leq\epsilon$.
Then,  for $t\in\mathbb{R}_+$ we get
{\small\begin{eqnarray}\nonumber
|({\mathcal{F}}x)(t)-({\mathcal{F}}y)(t)|\leq\left|\frac{f(t,x(t))}{\Gamma(\alpha)}\int_0^t
\frac{u(t,s,x(s),x(\lambda s))}{(t-s)^{1-\alpha}}\;ds\right.\;\;\;\;\;\;\;\;\;\;\;\;\;\;\;\;\;\;\;
\;\;\;\;\;\;\;\;\;\;\;\;\;\\
\nonumber\;\;\;-\left.\frac{f(t,y(t))}{\Gamma(\alpha)}
\int_0^t\frac{u(t,s,y(s),y(\lambda s))}{(t-s)^{1-\alpha}}\;ds\right|\\
\nonumber
\leq
\frac{|f(t,x(t))-f(t,y(t))|}{\Gamma(\alpha)}\int_0^t
\frac{|u(t,s,x(s),x(\lambda s))|}{(t-s)^{1-\alpha}}\;ds\;\;\;\;\;\;\;\;\;\;\;\;\\
\nonumber\;\;\;\;\;\;\;\;\;\;\;\;\;+\frac{|f(t,y(t))|}{\Gamma(\alpha)}
\int_0^t
\frac{|u(t,s,x(s),x(\lambda s))-u(t,s,y(s),y(\lambda s))|}{(t-s)^{1-\alpha}}ds\\
\nonumber
\leq
\frac{m(t)|x(t)-y(t)|}{\Gamma(\alpha)}\int_0^t
\frac{|u(t,s,x(s),x(\lambda s))-u(t,s,0,0)|+|u(t,s,0,0)|}{(t-s)^{1-\alpha}}\;ds\\
\nonumber\;\;\;\;\;\;\;\;\;+\frac{m(t)|y(t)|+|f(t,0)|}{\Gamma(\alpha)}
\int_0^t
\frac{n(t)\;\Phi(|x(s)-y(s)|,|x(\lambda s)-y(\lambda s)|)}{(t-s)^{1-\alpha}}ds\\
\nonumber
\leq
\frac{m(t)|x(t)-y(t)|}{\Gamma(\alpha)}\int_0^t
\frac{n(t)\;\Phi(|x(s)|,|x(\lambda s)|)+u^*(t)}{(t-s)^{1-\alpha}}\;ds\\
\nonumber\;\;\;\;\;\;\;\;\;+\frac{m(t)|y(t)|+|f(t,0)|}{\Gamma(\alpha)}
\int_0^t
\frac{n(t)\;\Phi(|x(s)|+|y(s)|,|x(\lambda s)|+|y(\lambda x)|)}{(t-s)^{1-\alpha}}ds\\
\nonumber
\leq
\frac{|x(t)-y(t)|}{\Gamma(\alpha+1)}[\phi(t)\Phi(\|x\|,\|x\|)+\psi(t)]\;\;\;\;\;\;\;\;\;\;\;\;\;\;\;\;\;
\;\;\;\;\;\;\;\;\;\;\;\;\;\;\;\;\;\;\\
\nonumber\;\;\;\;\;+\frac{\phi(t)|y(t)|{\Phi(\|x\|+\|y\|,\|x\|+\|y\|)}}{\Gamma(\alpha+1)}\;\;\;\;\;\;\;\;\;\;\;\;\;\;\;\;\;\;\\
\nonumber+
\frac{\xi(t){\Phi(\|x\|+\|y\|,\|x\|+\|y\|)}}{\Gamma(\alpha+1)}\;\;\;\;\;\;\;\;\;\;\;\;\;\;\;\;\;\;\;\;\;\;\;\;\;\;\;\\
\nonumber
\leq
\frac{1}{\Gamma(\alpha+1)}[\varepsilon\phi(t)\Phi(r_0,r_0)+\varepsilon\psi(t)\;\;\;\;\;
\;\;\;\;\;\;\;\;\;\;\;\;\;\;\;\;\;\;\;\;\;\;\;\;\;\;\;\;\;\;\;\;\;\;\;\;\;\;\;\;\;\;\;\\
\label{hh}\;\;\;\;
+\phi(t)r_0\Phi(2r_0,2r_0)+\xi(t)\Phi(2r_0,2r_0)].\;\;\;\;\;\;\;\;\;\;\;\;\;\;\;\;\;\;\;\;\;\;\;\;
\end{eqnarray}}
Thus $\mathcal{F}$ is continuous on $B_{r_0}$, thanks hypothesis $(h_4)$.

In what follows let us take a nonempty set $X\subset B_{r_0}$. Then, for arbitrary $x,\;y\in X$
and for a fixed $t\in\mathbb{R}_+$, calculating in the same way as in estimate $(\ref{hh})$, we obtain
{\small\begin{eqnarray} \nonumber
|({\mathcal{F}}x)(t)-({\mathcal{F}}y)(t)|&\leq&
\frac{|x(t)-y(t)|}{\Gamma(\alpha+1)}[\phi(t)\Phi(\|x\|,\|x\|)+\psi(t)]\\
\nonumber\;&\;&\;\;\;+\frac{{\Phi(\|x\|+\|y\|,\|x\|+\|y\|)}}{\Gamma(\alpha+1)}[\phi(t)\|y\|+\xi(t)]\\
\nonumber
&\leq&
\frac{\phi(t)\Phi(r_0,r_0)+\psi(t)}{\Gamma(\alpha+1)}|x(t)-y(t)|\\
\nonumber&\;&\;\;\;\;\;+\frac{\Phi(2r_0,2r_0)}{\Gamma(\alpha+1)}
[\phi(t)r_0+\xi(t)].
\end{eqnarray}}
Hence, we can easily deduce the following inequality
$${\rm diam}({\mathcal{F}}X)(t)\leq\frac{\phi(t)\Phi(r_0,r_0)+\psi(t)}{\Gamma(\alpha+1)}
{\rm diam}X(t)+\frac{\Phi(2r_0,2r_0)}{\Gamma(\alpha+1)}
[\phi(t)r_0+\xi(t)].$$
Now, taking into account hypothesis $(h_4)$ we obtain
\begin{equation}
\label{hm1}
c({\mathcal{F}}X)\leq k\;c(X),
\end{equation}
where we denoted $k=\frac{\phi^*\Phi(r_0,r_0)+\psi^*}{\Gamma(\alpha+1)}$. Obviously, in view of
hypothesis $(h_5)$ we have that $k<1$.

In what follows, let us take arbitrary numbers $\varepsilon>0$ and $T>0$. Choose a function $x\in X$ and
take $t_1,\;t_2\in[0,T]$ such that $|t_2-t_1|\leq\varepsilon$. Without loss of generality we can assume that $t_2>t_1$. Then,
taking into account our hypotheses and (\ref{es}), we have
{\small
\begin{eqnarray}\nonumber
|({\mathcal{F}}x)(t_2)-({\mathcal{F}}x)(t_1)|\leq|a(t_2)-a(t_1)|+
\left|(Fx)(t_2){({\mathcal{U}}x)(t_2)}-(Fx)(t_1)
{({\mathcal{U}}x)(t_2)}\right|\\
\nonumber\;\;\;\;\;\;\;\;\;\;\;\;\;\;
\;\;\;\;\;\;\;+
\left|(Fx)(t_1){({\mathcal{U}}x)(t_2)}-(Fx)(t_1)
{({\mathcal{U}}x)(t_1)}\right|\\
\nonumber\;\leq\omega^T(a,\varepsilon)+
\frac{|f(t_2,x(t_2))-f(t_1,x(t_1))|}{\Gamma(\alpha)}\;\;\;\;\;\;\;\;\;\;\;\;\;\;\;\;\;\;\;\;\;\;\;\;\;\;\;\\
\nonumber\times\int_0^{t_2}\frac{|u(t_2,s,x(s),x(\lambda s))-u(t_2,s,0,0)|+|u(t_2,s,0,0)|}
{(t_2-s)^{1-\alpha}}\;ds\\
\nonumber\;\;\;+\frac{|f(t_1,x(t_1))-f(t_1,0)|+|f(t_1,0)|}{\Gamma(\alpha+1)}\left\{
2\varepsilon^\alpha[\hat{n}(T)\Phi(\|x\|,\|x\|)+\hat{u}(T)]+T^\alpha\omega_{\|x\|}^T(u,\varepsilon)\right\}\\
\nonumber\;\leq\omega^T(a,\varepsilon)+
\frac{m(t_2)|x(t_2)-x(t_1)|+\omega_f^T(\varepsilon)}{\Gamma(\alpha)}\int_0^{t_2}\frac{n(t_2)\Phi(|x(s)|,|x(\lambda s)|)+u^*(t_2)}
{(t_2-s)^{1-\alpha}}\;ds\\
\nonumber\;\;\;+\frac{m(t_1)|x(t_1)|+|f(t_1,0)|}{\Gamma(\alpha+1)}\left\{
2\varepsilon^\alpha[\hat{n}(T)\Phi(\|x\|,\|x\|)+\hat{u}(T)]+T^\alpha\omega_{\|x\|}^T(u,\varepsilon)\right\}\\
\nonumber\leq\omega^T(a,\varepsilon)+\frac{t_2^\alpha}{\Gamma(\alpha+1)}[m(t_2)\omega^T(x,\varepsilon)+\omega_f^T(\varepsilon)]
[n(t_2)\Phi(r_0,r_0)+u^*(t_2)]\\
\nonumber\;\;\;\;\;\;+\frac{\hat{m}(T)r_0+\hat{f}(T)}{\Gamma(\alpha+1)}\left\{
2\varepsilon^\alpha[\hat{n}(T)\Phi(\|x\|,\|x\|)+\hat{u}(T)]+T^\alpha\omega_{\|x\|}^T(u,\varepsilon)\right\}\\
\nonumber\leq\omega^T(a,\varepsilon)+\frac{[\phi(t_2)\Phi(r_0,r_0)+\psi(t_2)]}{\Gamma(\alpha+1)}\omega^T(x,\varepsilon)
+\frac{T^\alpha\;\omega_f^T(\varepsilon)}{\Gamma(\alpha+1)}
[\hat{n}(T)\Phi(r_0,r_0)+\hat{u}(T)]\\
\nonumber\;\;\;\;\;\;+\frac{\hat{m}(T)r_0+\hat{f}(T)}{\Gamma(\alpha+1)}\left\{
2\varepsilon^\alpha[\hat{n}(T)\Phi(\|x\|,\|x\|)+\hat{u}(T)]+T^\alpha\omega_{\|x\|}^T(u,\varepsilon)\right\}\\
\nonumber\leq\omega^T(a,\varepsilon)+\frac{[\phi^*\Phi(r_0,r_0)+\psi^*]}{\Gamma(\alpha+1)}\omega^T(x,\varepsilon)
+\frac{T^\alpha\;\omega_f^T(\varepsilon)}{\Gamma(\alpha+1)}
[\hat{n}(T)\Phi(r_0,r_0)+\hat{u}(T)]\\
\label{es1}\;+\frac{\hat{m}(T)r_0+\hat{f}(T)}{\Gamma(\alpha+1)}\left\{
2\varepsilon^\alpha[\hat{n}(T)\Phi(\|x\|,\|x\|)+\hat{u}(T)]+T^\alpha\omega_{\|x\|}^T(u,\varepsilon)\right\},\;\;\;\;
\end{eqnarray}}
where we denoted
$$\omega_f^T(\varepsilon)=\sup\{|f(t_2,x(t_2))-f(t_1,x(t_1))|:t_1,\;t_2\in[0,T],\;|t_2-t_1|\leq\varepsilon,\;
x\in[-r_0,r_0]\},$$
$$\hat{m}(T)=\max\{m(t):t\in[0,T]\}$$
and
$$\hat{f}(T)=\max\{|f(t,0)|:t\in[0,T]\}.$$

Now, keeping in mind the uniform continuity of the function $f=f(t,x)$ on the set
$[0,T]\times[r_0,r_0]$ and the uniform continuity of the function $u=u(t,s,x,y)$ on the set
$[0,T]\times[0,T]\times[r_0,r_0]\times[r_0,r_0]$, from estimate $(\ref{es1})$ we derive the following one
$$\omega_0^T({\mathcal{F}}X)\leq k\;\omega_0^T(X).$$
Hence we have
\begin{equation}
\label{hm2}
\omega_0^\infty({\mathcal{F}}X)\leq k\;\omega_0^\infty(X).
\end{equation}
From (\ref{hm1}) and (\ref{hm2}) and the definition of the
measure of noncompactness $\mu$ given by formula (\ref{f1}), we obtain
\begin{equation}
\label{nh}\mu({\mathcal{F}}X)\leq k\;\mu(X).
\end{equation}
In the sequel let us put $B_{r_0}^1={\rm Conv}{\mathcal{F}}(B_{r_0})$, $B_{r_0}^2={\rm Conv}{\mathcal{F}}(B_{r_0}^1)$
and so on. In this way we have constructed a decreasing sequence of nonempty, bounded, closed and convex subsets $(B_{r_0}^n)$ of
$B_{r_0}$ such that ${\mathcal{F}}(B_{r_0}^{n})\subset B_{r_0}^n$ for $n=1,\;2,\;\ldots$. Since the above reasons leading to
(\ref{nh}) holds for any subset X of $B_{r_0}$ we have
$$\mu(B_{r_0}^n)\leq k^n\mu(B_{r_0}),\;{\rm for}\;{\rm any}\;n=1,\;2,\;3,\;\ldots.$$
This implies that $\lim\limits_{n\rightarrow\infty}\mu(B_{r_0}^n)=0$. Hence, taking into account Definition\ref{d1} we infer
that the set $Y=\bigcap\limits_{n=1}^{\infty}B_{r_0}^n$ is nonempty, bounded, closed and convex subset of $B_{r_0}$. Moreover,
$Y\in{\rm ker}\mu$. Also, the operator $\mathcal{F}$ maps $Y$ into itself.

We will prove that the operator $\mathcal{F}$ is continuous on the set $Y$. In order to do this let us fix a number
$\varepsilon>0$ and take arbitrary functions $x,\;y\in Y$ such that $\|x-y\|\leq\varepsilon$. Keeping in mind the facts that
$Y\in{\rm ker}\mu$ and the structure of sets belong to ${\rm ker}\mu$ we can find a number $T>0$ such that for each $z\in Y$ and
$t\geq T$ we have that $|z(t)|\leq\varepsilon$. Since $\mathcal{F}$ maps $Y$ into itself we have that
${\mathcal{F}}x,\;{\mathcal{F}}y\in Y$. Thus, for $t\geq T$ we get
\begin{equation}
\label{s1}
|({\mathcal{F}}x)(t)-({\mathcal{F}}y)(t)|\leq|({\mathcal{F}}x)(t)|+|({\mathcal{F}}y)(t)|\leq 2\varepsilon.
\end{equation}
On the other hand, let us assume $t\in[0,T]$. Then we obtain

{\small\begin{eqnarray}\nonumber
|({\mathcal{F}}x)(t)-({\mathcal{F}}y)(t)|\leq\frac{m(t)|x(t)-y(t)|}{\Gamma(\alpha)}\int_0^t
\frac{n(t)\;\Phi(|x(s)|,|x(\lambda s)|)+u^*(t)}{(t-s)^{1-\alpha}}\;ds\\
\nonumber\;\;\;\;\;\;\;\;\;+\frac{m(t)|y(t)|+|f(t,0)|}{\Gamma(\alpha)}
\int_0^t
\frac{n(t)\;\Phi(|x(s)-y(s)|,|x(\lambda s)-y(\lambda s)|)}{(t-s)^{1-\alpha}}ds\\
\nonumber
\leq
\frac{[m(t)n(t)\Phi(r_0,r_0)+m(t)u^*(t)]\varepsilon}{\Gamma(\alpha)}\int_0^t
\frac{ds}{(t-s)^{1-\alpha}}\\
\nonumber\;\;\;\;\;\;\;\;\;+\frac{[m(t)n(t)r_0+n(t)|f(t,0)]\Phi(\varepsilon,\varepsilon)}{\Gamma(\alpha)}
\int_0^t
\frac{ds}{(t-s)^{1-\alpha}}\\
\nonumber
\leq
\frac{\phi(t)\Phi(r_0,r_0)+\psi(t)}{\Gamma(\alpha+1)}\;\varepsilon+\frac{\phi(t)r_0+\xi(t)}
{\Gamma(\alpha+1)}\Phi(\varepsilon,\varepsilon)\\
\label{s2}
\leq
\frac{\phi^*\Phi(r_0,r_0)+\psi^*}{\Gamma(\alpha+1)}\;\varepsilon+\frac{\phi^*r_0+\xi^*}
{\Gamma(\alpha+1)}\Phi(\varepsilon,\varepsilon).\;\;\;\;\;\;\;\;\;
\end{eqnarray}}
Now, taking into account $(\ref{s1})$, $(\ref{s2})$ and hypothesis $(h_4)$ we conclude that the operator
$\mathcal{F}$ is continuous on the set $Y$.

Finally, linking all above obtained facts about the set $Y$ and the operator ${\mathcal{F}}:Y\rightarrow Y$ and using
the classical Schauder fixed point principal we deduce that the operator $\mathcal{F}$ has at least one fixed point
$x$ in the set $Y$. Obviously the function $x=x(t)$ is a solution of the quadratic integral equation $(\ref{e1})$. Moreover,
since $Y\in{\rm ker}\mu$ we have that $x(t)\rightarrow 0$ as $t\rightarrow\infty$. This completes the proof.
\section{EXAMPLE}
Consider the following quadratic integral equation of fractional order with linear
modification of the argument, $\alpha=\frac{1}{2}$,
\begin{equation}\label{ex2}
x(t)=te^{-t}+\frac{t+t^2x(t)}{2\Gamma(\frac{1}{2})}\int_0^t
\frac{(|x(s)|+|x(\lambda s)|)e^{-2t-s}+1\slash(1+5t^{5\slash 2})}{\sqrt{t-s}}\;ds.
\end{equation}

In this example, we have that $a(t)=te^{-t}$ and this function satisfies
hypothesis $(h_1)$ and $\|a\|=a(1)=1\slash e$. Moreover, $f(t,x(t))=(t+t^2x(t))\slash 2$ and
satisfies hypothesis $(h_2)$ with $m(t)=t^2\slash 2$,
and $|f(t,0)|=f(t,0)=t\slash 2$. Also, $u(t,s,x,y)=(x+y)e^{-2t-s}+1\slash(1+5t^{5\slash 2})$
verifies hypothesis $(h_3)$ with $n(t)=e^{-2t}$, $\phi(x,y)=x+y$ and
$u(t,s,0,0)=1\slash(1+5t^{5\slash 2})$.
Now, we have
$$\phi(t)=\frac{1}{2}t^{5\slash 2}e^{-2t},$$
$$\psi(t)=\frac{t^{5\slash 2}}{2(1+5t^{5\slash 2})},$$
$$\xi(t)=\frac{1}{2}t^{3\slash 2}e^{-2t}$$
and
$$\eta(t)=\frac{t^{3\slash 2}}{2(1+5t^{5\slash 2})}.$$
It is easy to see that the functions $\phi,\;\psi,\;\xi$ and $\eta$ are bounded on $\mathbb{R}_+$ and
also $\lim\limits_{t\rightarrow\infty}\phi(t)=\lim\limits_{t\rightarrow\infty}\xi(t)=0.$ Hence,
hypothesis $(h_4)$ is satisfied. Moreover, we have
$$\phi^*=\phi(5\slash 4)=\frac{1}{2}(5\slash 4)^{5\slash 2}e^{-5\slash 2}=0.0716982...,$$
$$\psi^*=0.1,$$
$$\xi^*=\phi(3\slash 4)=\frac{1}{2}(3\slash 4)^{5\slash 2}e^{-3\slash 2}=0.0543477...$$
and
$$\eta^*=\eta((0.2)^{2\slash 5})=0.0410503...\;.$$
In this case the inequality $(\ref{inq})$ has the form
\begin{equation}
\label{iem}
\Gamma(\frac{3}{2})e^{-1}+2r^2\phi^*+r\psi^*+2r\xi^*+\eta^*\leq r\Gamma(\frac{3}{2}).
\end{equation}
Let us denote by $H(r)$ the left hand side of the last inequality, i.e.
$$H(r)=\Gamma(\frac{3}{2})e^{-1}+2r^2\phi^*+r\psi^*+2r\xi^*+\eta^*.$$
For $r=1$ we obtain
\begin{eqnarray}
\nonumber
H(1)&=&\Gamma(\frac{3}{2})e^{-1}+2\phi^*+\psi^*+2\xi^*+\eta^*\\
\nonumber&=&\Gamma(\frac{3}{2})\;0.3678794...+0.3931423...\;.
\end{eqnarray}
Hence, inequality (\ref{iem}) admits $r_0=1$ as a
positive solution since $\Gamma(\frac{3}{2})\simeq 0.886227$.
Moreover,
$$\phi^*\phi(r_0,r_0)+\psi^*\simeq0.2433966<\Gamma(\frac{3}{2}).$$
Therefore, Theorem $\ref{mth}$ guarantees that equation (\ref{ex2}) has a
solution $x=x(t)$ in the space $\mathbb{R}_+$ belonging to the ball $B_1$ such that
$x(t)\rightarrow0$ as
$t\rightarrow\infty$.

\end{document}